\documentclass[11pt,bezier]{article}
\usepackage{amsmath,graphicx,amssymb,amsfonts,color}

\font\bg=cmbx10 scaled\magstep1 \font\Bg=cmbx12 scaled\magstep3
\font\small=cmr8

\newtheorem{newlemma}{{\bf Lemma}}

\newenvironment{lema}{\begin{newlemma}{\hspace{-0.5
em}{\bf.}}}{\end{newlemma}}

\newtheorem{newteorem}{{\bf Theorem}}

\newenvironment{teorem}{\begin{newteorem}{\hspace{-0.5
em}{\bf.}}}{\end{newteorem}}

\newtheorem{newkorolari}{{\bf Corollary}}

\newenvironment{korolari}{\begin{newkorolari}{\hspace{-0.5
em}{\bf.}}}{\end{newkorolari}}

\newtheorem{newkonjek}{{\bf Conjecture}}

\newenvironment{konjek}{\begin{newkonjek}{\hspace{-0.5
em}{\bf.}}}{\end{newkonjek}}

\newtheorem{newproblem}{{\bf Problem}}

\newenvironment{problem}{\begin{newproblem}{\hspace{-0.5
em}{\bf.}}}{\end{newproblem}}

\begin{document}
\tolerance=10000 \baselineskip18truept
\newbox\thebox
\global\setbox\thebox=\vbox to 0.2truecm{\hsize
0.15truecm\noindent\hfill}
\def\boxit#1{\vbox{\hrule\hbox{\vrule\kern0pt
     \vbox{\kern0pt#1\kern0pt}\kern0pt\vrule}\hrule}}
\def\qed{\lower0.1cm\hbox{\noindent \boxit{\copy\thebox}}\bigskip}
\def\ss{\smallskip}
\def\ms{\medskip}
\def\nt{\noindent}

\centerline{\Bg Graphs which their certain polynomials  }
 \vspace{.3cm}

\centerline {\Bg  have few distinct roots- a survey}
\bigskip
\baselineskip12truept \centerline{\bg Saeid
Alikhani$^{a,b,}$\footnote{\baselineskip12truept\it\small  E-mail:
alikhani@yazduni.ac.ir}} 
\baselineskip20truept \centerline{\it $^a$Department of
Mathematics, Yazd University} \vskip-8truept \centerline{\it
89195-741, Yazd, Iran}
 \vskip-5truept \centerline{\it
$^{b}$School of Mathematics, Institute for Research in Fundamental
Sciences (IPM)} \vskip-8truept \centerline{\it P.O. Box
19395-5746, Tehran, Iran}

\nt\rule{15cm}{0.1mm}

\nt{\bg ABSTRACT}
 \medskip

 \nt
Let $G=(V,E)$ be a simple graph. We consider domination
polynomial, matching  polynomial and edge cover polynomial of $G$.
Graphs which their polynomials have few roots can give sometimes a
very surprising information about the structure of the graph. In
this paper we study  graphs which their domination polynomial,
independence polynomial and edge cover polynomial have few roots.

\nt {\bf Keywords:} Domination polynomial, Edge cover polynomial;
Matching polynomial; Independence polynomial;   Root.

\nt{\bf Msc:} 05C69, 11B83.                        

\nt\rule{15cm}{0.1mm}

 \section{Introduction}

\nt Let $G=(V,E)$ be a simple graph.    Graph polynomials are a
well-developed area useful for analyzing properties of graphs. The
study of  graphs which their polynomials have few roots  can give
sometimes a very surprising information about the structure of the
graph.

\nt We consider domination polynomial, matching polynomial (and
independence polynomial) and edge cover polynomial of a graph $G$.
  For convenience, the
definition of these polynomials  will be given in the next
sections.

\nt The {\it corona} of two graphs $G_1$ and $G_2$, as defined by
Frucht and Harary in~\cite{harary},
 is the graph
$G=G_1 \circ G_2$ formed from one copy of $G_1$ and $|V(G_1)|$
copies of $G_2$, where the ith vertex of $G_1$ is adjacent to
every vertex in the ith copy of $G_2$. The corona $G\circ K_1$, in
particular, is the graph constructed from a copy of $G$,
 where for each vertex $v\in V(G)$, a new vertex $v'$ and a pendant edge $vv'$ are added.
 The {\it join} of two graphs $G_1$ and $G_2$, denoted by $G_1\vee G_2$ is a graph with vertex set  $V(G_1)\cup V(G_2)$
 and edge set $E(G_1)\cup E(G_2)\cup \{uv| u\in V(G_1)$ and $v\in V(G_2)\}$.

\nt The characterization of graphs with few distinct roots of
characteristic polynomials (i.e. graphs with few distinct
eigenvalues) have been the subject of many researches. Graphs with
three adjacency eigenvalues have been studied by Bridges and Mena
\cite{bridges} and Klin and Muzychuk \cite{muz}. Also van Dam
studied graphs with three and four distinct eigenvalues
\cite{dam1,dam2,dam3,dam4,dam5}. Graphs  with three distinct
eigenvalues and index less than $8$ studied by Chuang and Omidi in
\cite{omidi}.

\nt In this paper we study graphs which their domination
polynomial, matching polynomial (and so independence polynomial)
and edge cover polynomial have few roots. In Section 2 we
investigate graphs with few domination roots. In Section 3 we
characterize graphs which their matching polynomials (and so their
independence polynomials) have few roots.  Finally in Section 4 we
characterize graphs with few edge cover roots.

\nt As usual we use $\lfloor x\rfloor,\lceil x \rceil$ for the
largest integer less than or
  equal to $x$ and for the smallest integer greater than or equal to $x$, respectively.
In this paper we denote the set $\{1,2,\ldots,n\}$ simply by $[n]$

\section{Graphs with few domination roots}

\nt Let $G=(V,E)$ be a graph of order $|V|=n$. For any vertex
$v\in V$, the {\it open neighborhood} of $v$ is the set $N(v)=\{u
\in V|uv\in E\}$ and the {\it closed neighborhood} of $v$ is the
set $N[v]=N(v)\cup \{v\}$.
 For a set $S\subseteq V$, the open neighborhood of $S$
 is $N(S)=\bigcup_{v\in S} N(v)$ and the closed neighborhood of $S$ is $N[S]=N(S)\cup S$.
A set $S\subseteq V$ is a dominating set if $N[S]=V$, or
equivalently,
 every vertex in $V\backslash S$ is adjacent to at least one vertex in $S$.
 The domination number $\gamma(G)$ is the minimum cardinality of a dominating set in $G$.
For a detailed treatment of this parameter, the reader is referred
to~\cite{domination}. Let ${\cal D}(G,i)$ be the family of
dominating sets of a graph $G$ with cardinality $i$ and let
$d(G,i)=|{\cal D}(G,i)|$.
 The {\it domination polynomial} $D(G,x)$ of $G$ is defined as
$D(G,x)=\sum_{i=\gamma(G)}^{|V(G)|} d(G,i) x^{i}$, where
$\gamma(G)$ is the domination number of $G$
(\cite{saeid1,saeid3}). The path $P_{4}$ on 4 vertices, for
example, has one dominating  set of cardinality 4, four dominating
 sets of cardinality 3, and four dominating  sets of cardinality
2; its domination polynomial is then $D(P_{4},x)=x^4+4x^3+4x^2$.

 \nt A root of $D(G,x)$ is called a {\it domination  root} of $G$. In this section, the set
of distinct roots of $D(G,x)$ is denoted by $Z(D(G,x))$.

\ms

\nt We recall  a formula for the computation of the domination
polynomial of  join of two graphs.

\begin{teorem}\label{theorem8}{\rm(\cite{saeid2})}
Let $G_1$ and $G_2$ be  graphs of orders $n_1$ and $n_2$,
respectively. Then
\[
D(G_1\vee
G_2,x)=\Big((1+x)^{n_1}-1\Big)\Big((1+x)^{n_2}-1\Big)+D(G_1,x)+D(G_2,x).
\]
\end{teorem}

\nt The following theorem is an easy result about roots of $K_n$
and $K_{1,n}$.

\begin{teorem}\label{theorem6.2.1}{\rm(\cite{zero})}
\begin{enumerate}
\item[(i)]
 For every $n\in {\mathbb N}$,
 \[
 Z\Big(D(K_n,x)\Big)=\Big\{ \exp({\frac{2k\pi i}{n}})-1 \hspace{.4mm} \Big |\hspace{.4mm}  k=0,1,\ldots,n-1\Big\}.
 \]
\item[(ii)]
 For every $n\in {\mathbb N}$, $D(K_{1,n},x)$ has exactly two real roots for odd $n$ and exactly three real roots for even $n$.
\end{enumerate}

\end{teorem}

\nt The following result is about the domination roots of
$K_{n,n}$:

\begin{teorem}\label{theorem6.2.2}
For every even $n$, no nonzero real numbers is  domination root of
$K_{n,n}$.
\end{teorem}
\nt{\bf Proof.}  By Theorem~\ref{theorem8} for
$G_1=G_2=\overline{K_n}$, we have
\[
D(K_{n,n},x)=\Big((1+x)^n-1\Big)^2+2x^n.
\]
If $D(K_{n,n},x)=0$, then $\Big((1+x)^n-1\Big)^2=-2x^n$.
 Obviously this equation does not have real nonzero solution for even $n$.\quad\qed

\ms

\nt In~\cite{zero,saeid1} we characterized graphs with one, two
and three distinct   domination roots. Since $0$ is a root of any
domination polynomial of graph $G$, we have the following theorem.

\begin{teorem} \label{theorem6'}{\rm(\cite{saeid2})}
A graph $G$ has one domination root if and only if $G$ is a union
of isolated vertices.
\end{teorem}

\nt The following theorem characterize graphs with two distinct
domination roots.

\begin{teorem}\label{theorem6} {\rm(\cite{saeid1})}
Let $G$ be a connected graph with exactly two distinct domination
roots. Then there exists natural number $n$ such that
$D(G,x)=x^n(x+2)^n$. Indeed $G=H\circ K_1$ for some graph  $H$ of
order $n$. Moreover, for every graph $H$ of order $n$, $D(H\circ
K_1,x)=x^n(x+2)^n$.
\end{teorem}


\begin{teorem}\label{theorem5'} {\rm(\cite{saeid1})}
Let $G$ be a connected graph of order $n$. Then, $Z(D(G,x))= \{0,
\frac{-3\pm \sqrt{5}}{2}\}$, if and only if $G=H\circ
\overline{K}_2$, for some graph $H$. Indeed
$D(G,x)=x^{\frac{n}{3}}(x^2+3x+1)^{\frac{n}{3}}$.
\end{teorem}

\nt Theorem \ref{theorem5'} characterize graphs with $Z(D(G,x))=
\{0, \frac{-3\pm \sqrt{5}}{2}\}$ (see Figure \ref{figure1}). The
following  theorem shows that roots of graphs with exactly three
distinct domination roots can not be any numbers.

\begin{teorem}\label{theorem10} {\rm(\cite{saeid1})}
For every graph $G$ with exactly three distinct domination roots
\[
Z(D(G,x))\subseteq  \Big\{0,\frac{-3\pm\sqrt{5}}{2}, -2\pm
\sqrt{2}i ,\frac{-3\pm\sqrt{3}i}{2}\Big\}.
\]
\end{teorem}

\nt {\bf Remark.} Since
$Z(D(K_3,x))=\{0,\frac{-3\pm\sqrt{3}i}{2}\}$,
$Z(D(P_3,x))=\{0,\frac{-3\pm\sqrt{5}}{2}\}$ and
$Z(D(C_4,x))=\{0,-2\pm \sqrt{2}i\}$, for every number of set
$\Big\{0,\frac{-3\pm\sqrt{5}}{2}, -2\pm \sqrt{2}i
,\frac{-3\pm\sqrt{3}i}{2}\Big\}$, there exist a graph which its
domination polynomial have exactly three distinct roots. Since
cycles are determined by their domination polynomials
(\cite{arsob}), so graphs with exactly three domination roots from
$\{ 0,\frac{-3\pm\sqrt{5}}{2}, -2\pm \sqrt{2}i \}$ are $C_4$ and
graphs of the form $H\circ \overline{K}_2$, for some graph $H$.


\nt By Theorems \ref{theorem6'},\ref{theorem6} and \ref{theorem10}
we have the following corollary:

\begin{korolari}
For every graph $G$ with at most three distinct domination roots
\[
Z(D(G,x))\subseteq  \Big\{0,-2,\frac{-3\pm\sqrt{5}}{2}, -2\pm
\sqrt{2}i ,\frac{-3\pm\sqrt{3}i}{2}\Big\}.
\]
\end{korolari}

 \nt Now we shall study graphs with exactly four distinct
 domination roots.

\nt Let $G_n$ be an arbitrary graph of order $n$. Let to denote
the graph $G_n\circ K_1$ simply by $G_n^*$.  Here we consider the
labeled $G_n^*$ as show in Figure~\ref{figure2} (the graph in this
figure is $P_n^*$ which called centipede). We denote the graph
obtained from $G_n^*$ by deleting the vertex labeled $2n$ as
$G_n^*-\{2n\}$.

\ms



\nt The following theorem state  a  recursive formula for the
domination polynomial of $G_n^*-\{2n\}$.

\begin{teorem}\label{theorem50} {\rm(\cite{four})}
 For every $n\geq 5$,
\[
D(G_{n}^*-\{2n\},x)=x\Big[D(G_{n-1}^*,x)+D(G_{n-2}^*,x)\Big]+x^2D(G_{n-2}^*,x).
\]
\end{teorem}

\nt The following theorem give the formula for
$D(G_n^*-\{2n\},x)$.
\begin{teorem}\label{theorem7} {\rm(\cite{four})}
 For every $n\geq2$, $D(G_n^*-\{2n\},x)=(x^2+3x+1)x^{n-1}(x+2)^{n-2}$.
\end{teorem}

\ms \nt The following theorem characterize graphs with  four
domination roots $-2,0, \frac{-3\pm \sqrt{5}}{2}$.

\begin{teorem} {\rm(\cite{four})}
Let $G$ be a connected graph of order $n$. Then, $Z(D(G,x))=
\{-2,0, \frac{-3\pm \sqrt{5}}{2}\}$, if and only if
$G=G_{\frac{n}{2}}^*-\{n\}$, for some graph $G_{\frac{n}{2}}$ of
order $\frac{n}{2}$. Indeed
\[
D(G,x)=(x^2+3x+1)x^{\frac{n}{2}-1}(x+2)^{\frac{n}{2}-2}.
\]
\end{teorem}


\ms


\nt Using tables of domination polynomials (see \cite{thesis}), we
think that numbers which are roots of graphs with exactly four
distinct domination roos are finite and are about nine numbers,
but we are not able to prove it. So complete characterization of
graphs with exactly four distinct domination roots remains as open
problem.

\nt It is natural to ask about the domination roots of paths and
cycles. First we recall the following theorem:

\begin{teorem}\label{theorem3}
For every $n\geq 3$,

\begin{enumerate}
\item[(i)] {\rm(\cite{saeid4})}
$D(P_n,x)=x\Big[D(P_{n-1},x)+D(P_{n-2},x)+D(P_{n-3},x)\Big],$

\item[(ii)]{\rm(\cite{saeid3})}
$D(C_n,x)=x\Big[D(C_{n-1},x)+D(C_{n-2},x)+D(C_{n-3},x)\Big]$,
\end{enumerate}
\end{teorem}

\nt We think that real roots of the families $D(P_{n},x)$ and
$D(C_{n},x)$ are dense in the interval $[-2,0]$, for $n\geq 4$,
but we couldn't to prove it yet.

\section{Graphs with few matching roots}

\nt In this section we study graphs which their matching
polynomials have few roots. First we state  the definition of
matching polynomial.  Let $G =(V,E)$ be a graph of order $n$ and
size $m$. An $r$-matching of $G$ is a set of $r$ edges of $G$
which no two of them have common vertex. The maximum number of
edges in a matching of a graph $G$ is called the matching number
of $G$ and denoted by
$\alpha'(G)$. 
The matching polynomial  is defined by
$$\mu(G,x)=\sum_{k=0}^{\lfloor\frac{n}{2}\rfloor} (-1)^k m(G,k) x^{n-2k},$$
where $m(G,k)$ is the number of $k$-matching of $G$ and
$m(G,0)=1$. The roots of $\mu(G,x)$ are called the matching roots
of $G$. As an example the matching polynomial of path  $P_5$ is
$\mu(P_5, x) = (x-1)x(x+1)(x^2-3)$. For more details of this
polynomial refer to \cite{match1,match2,match3}.

\nt Two  following theorems may are the first results on matching
roots.

\begin{teorem} {\rm(\cite{hel})}
 The roots of matching polynomial of any graph are all real
numbers.

\end{teorem}

\begin{teorem} {\rm(\cite{godsil})}
If $G$ has a Hamiltonian path, then all roots of its matching
polynomial are simple (have multiplicity 1).
\end{teorem}

\ms

\nt We need the following definition to study graphs with few
matching roots.

\ms \nt
Add a single vertex $u$ to the graph $rK_{1,k}\cup tK_1$ and join
$u$ to the other vertices by $p+q$ edges so that the resulting
graph is connected and $u$ is adjacent with $q$ centers of the
stars (for $K_{1,1}$ either of the vertices may be considered as
center). We denote the resulting graph by $\mathcal{S}(r, k, t; p,
q)$ (see Figure \ref{figure3}). Clearly $r + t \leq p + q \leq r(k
+ 1) + t$ and $0\leq q \leq r$ (see \cite{ghorbani}).


\nt For any $G\in \mathcal{S}(r, 3, t; p, q)$, we add $s$ copies
of $K_3$ to $G$ and join them by $l$ edges to the vertex $u$ of
$G$. Clearly $s\leq l \leq 3s$. We denote the set of these graphs
by $\mathcal{H}(r, s, t; p, q,l)$.

\nt The following theorem gives the matching polynomial of graph
$G$ in the family $\mathcal{S}(r, k, t; p, q)$.

\begin{teorem}{\rm(\cite{ghorbani})}
For every $G\in \mathcal{S}(r, k, t; p, q)$, $$\mu(G, x)=
x^{r(k-1)+t-1}(x^2- k)^{r-1}( x^4-(p + k)x^2 + (p- q)(k-1) + t).$$
\end{teorem}

\begin{teorem}{\rm(\cite{ghorbani})}
For every $G\in \mathcal{H}(r, s, t; p, q,l)$, $$\mu(G, x)=
x^{2r+s+t-1}(x^2- 3)^{r+s-1}( x^4-(p +l+3)x^2 +3t+2 (p-t-q)+l).$$
\end{teorem}

\nt Similar to \cite{ghorbani} we distinguish some special graphs
in the families $S$ and $H$ which are important for our study.

\nt We denote the family $\mathcal{S}(r, 1, 0; s, q)$ which
consists  a single graph by $S(r, s)$. Note that in this case $q$
is determined by $r$ and $s$, namely $q = s -r$. Its matching
polynomial is
$$\mu(S(r, s), x) = x(x^2 - s- 1)(x^2- 1)^{r-1}.$$

\nt The family $\mathcal{S}(r,k, 0; r, r)$ consists of a single
graph which is denoted by $T(r, k)$. Its matching polynomial is
$$\mu(T(r, k), x) = x^{r(k-1)+1}(x^2-r-k)(x^2-k)^{r-1}.$$

\nt We also denote the unique graphs in $\mathcal{S}(1, k, t; l +
t, 0)$ and $\mathcal{S}(1, k,t; l+ t + 1, 1)$ by $K(k, t;l)$ and
$K'(k, t;l)$, respectively. Their matching polynomials are
$$\mu(K(k, t;l), x) = x^{k+t-2}(x^4-(k + t +l)x^2 + (l+ t)(k-1) +
t);$$
$$\mu(K'(k,t;l), x) = x^{k+t-2}(x^4-(k + t +l+ 1)x^2 + (l + t)(k- 1) + t
).$$

\nt Moreover, we denote the unique graph in $\mathcal{H}(0, 1, t;
t, 0,l)$ by $L(t,l)$ for $l= 1,2,3$. We have $\mu(L(t,l), x) =
x^t(x^4-(t +l + 3)x^2 + 3t +l)$.  Typical graphs from the above
families are shown in Figure \ref{figure3'}.


\begin{teorem} {\rm(\cite{ghorbani})}
 Let $G$ be a connected graph and $z(G)$ be the number
of its distinct matching roots.
\begin{enumerate}

\item[(i)] If $z(G) = 2$, then $G\simeq K_2$.

\item[(ii)] If $z(G) = 3$, then $G$ is either a star or $K_3$.

\item[(iii)] If $z(G) = 4$, then $G$ is a non-star graph with $4$
vertices.

\item[(iv)] If $z(G) = 5$, then $G$ is one of the graphs
$K(k,t;l), K'(k,t;l), L(t;l), T(r,k), S(r, s)$, for some integers
$k,r, s, t, l$ or a connected non-star graph with $5$ vertices.
\end{enumerate}

\end{teorem}

\nt Using above theorem  we would like to study graphs with few
independence roots. First we recall the definition of independence
polynomial.

\ms \nt
 An independent set of a graph $G$ is a
set of vertices where no two vertices are adjacent. The
independence number is the size of a maximum independent set in
the graph and denoted by $\alpha(G)$. For a graph $G$, let $i_{k}$
denote the number of independent sets of cardinality $k$ in $G$
($k=0,1,\ldots,\alpha$). The independence polynomial of $G$,
\[
I(G,x)=\sum_{k=0}^{\alpha} i_{k} x^{k},
\]
is the generating polynomial for the independent sequence
$(i_{0},i_{1},i_{2},\ldots,i_{\alpha})$. For more study on
independence polynomial and independence root refer to
\cite{jamc,ind1,ind2}.

\nt The path $P_{4}$ on 4 vertices, for example, has one
independent set of cardinality 0 (the empty set), four independent
sets of cardinality 1, and three independent sets of cardinality
2; its independence polynomial is then $I(P_{4},x)=1+4x+3x^{2}$.

\ms

\nt Here we recall the definition of line graph.
 Given a graph $H = (V, E)$, the line graph of $H$,
denoted by $L(H)$, is a graph with vertex set $E$, two vertices of
$L(H)$ are adjacent if and only if the corresponding edges in $H$
share at least one endpoint. We say that $G$ is a line graph if
there is a graph $H$ for which $G = L(H)$.

\begin{teorem} \label{theorem11} {\rm(\cite{godsil2})}
For every graph $G$, $\mu(G,x)=x^n I(L(G),-\frac{1}{x^2})$.
\end{teorem}

\nt The following corollary is an immediate consequence of Theorem
\ref{theorem11}.

\begin{korolari} \label{corollary2}
If $\alpha\neq 0$ is a  matching root of $G$, then
$-\frac{1}{\alpha^2}$ is an independence  root of $L(G)$.
\end{korolari}

\nt Now we are ready to state a theorem for graphs which its
independence polynomial have  few roots. The following theorem
follows from Theorem \ref{theorem11} and Corollary
\ref{corollary2}.

\begin{teorem}\label{theorem18}
 Let $G$ be a connected graph and $z(L(G))$ be the number
of  distinct independence  non-zero roots of $L(G)$.
\begin{enumerate}

\item[(i)] If $z(L(G)) = 2$, then $L(G)\simeq K_2$.

\item[(ii)] If $z(L(G)) = 3$, then $L(G)$ is either a star or
$K_3$.

\item[(iii)] If $z(L(G)) = 4$, then $L(G)$ is a non-star graph
with $4$ vertices.

\item[(iv)] If $z(L(G)) = 5$, then $L(G)$ is one of the graphs
then $G$ is one of the graphs $K(k,t;l), K'(k,t;l), L(t;l),
T(r,k), S(r, s)$, for some integers $k,r, s, t, l$ or a connected
non-star graph with $5$ vertices.
\end{enumerate}

\end{teorem}

\begin{teorem} \label{theorem5}{\rm(\cite{jamc})}
\begin{enumerate}
\item[(i)]
 For any integer $n$, $I(P_{n},x)$ has the following
zeros,
\[
p_{s}^{(n)}=-\frac{1}{2\left(1+{\displaystyle
\cos\frac{2s\pi}{n+2}}\right)}, \quad s=1,2,\ldots,\lfloor
\frac{n+1}{2}\rfloor.
\]

\item[(ii)] For any integer $n\geq3$, $I(C_{n},x)$ has the
following zeros,
\[
c_{s}^{(n)}=-\frac{1}{2\left(1+{\displaystyle
\cos\frac{(2s-1)\pi}{n}}\right)}, \quad s=1,2,\ldots,\lfloor
\frac{n}{2}\rfloor.
\]
\end{enumerate}
\end{teorem}

\begin{korolari} \label{corollary4}{\rm(\cite{jamc})}
The independence roots of the family $\{P_n\}$ and $\{C_{n}\}$ are
real and dense in $(-\infty,-\frac{1}{4}]$.
\end{korolari}

\nt In \cite{rational} authors studied  graphs whose independence
roots are rational.

\begin{teorem}{\rm(\cite{rational})}
Let $G$ be a graph with rational polynomial $I(G, x)$. If
$i(G,r-1)> i(G, r)$, for some $r$, $1\leq r \leq \alpha(G)$, then
$G$ has an independence root in the interval $(-1,
\frac{-1}{\alpha(G)}]$.
\end{teorem}

\nt The following theorem characterize graphs with exactly one
independence roots:

\begin{teorem}{\rm(\cite{rational})}
Let $G$ be a graph of order $n$. Then $I(G,x)$ has exactly one
root if and only if $G=rK_s$, where $n=rs$ for some natural $r,
s$.
\end{teorem}

\section{Graphs with few edge-cover roots}

\nt In this section we characterize  graphs which their edge cover
polynomials have one and two distinct roots. First we state the
definition of edge-cover polynomial of a graph.

\ms
 \nt For every graph $G$ with no isolated vertex, an edge covering of
$G$ is a set of edges of $G$ such that every vertex is incident to
at least one edge of the set. A minimum edge covering is an edge
cover of the smallest possible size. The edge covering number of
$G$ is the size of a minimum edge cover of $G$ and denoted by
$\rho(G)$.  The edge cover polynomial of $G$ is the polynomial
$E(G,x)=\sum_{i=1}^m e(G,i)x^i$, where $e(G,i)$ is the number of
edge covering sets of $G$ of size $i$. Note that if graph $G$ has
isolated vertex then we put $E(G,x)=0$ and if
$V(G)=E(G)=\emptyset$, then $E(G,x)=1$. For more detail on this
polynomial refer to \cite{oboudi1,oboudi2}.

\nt As an example the edge cover polynomial of path  $P_5$ is
$E(P_5, x) =2x^3+x^4$. Also $E(K_{1,n},x)=x^n$.

\ms

\nt The following results are about edge cover polynomial of a
graph:

\begin{lema}\label{lemma2} {\rm(\cite{oboudi1})}
 Let $G$ be a graph of order $n$ and size $m$ with no
isolated vertex. If the edge cover polynomial of $G$ is
$E(G,x)=\sum_{i=\rho(G)}^m  e(G,i)x^i$, then the following hold:

\begin{enumerate}
\item[(i)] $n\leq 2 \rho(G)$. \item[(ii)] If $i_0 = min
\{i|e(G,i)={m \choose i}\}$, then $\delta = m- i_0 + 1$.

\item[(iii)] If $G$ has no connected component isomorphic to
$K_2$, then $a_{\delta(G)} ={m \choose m-\delta}- e(G,m -\delta)$.
\end{enumerate}
\end{lema}

\begin{teorem}{\rm(\cite{oboudi2})}
Let $G$ be a graph. Then $E(G,x)$ has at least $\delta(G)-2$
non-real roots (not necessary distinct). In particular, when
$\delta(G) = 3$, $E(G,x)$ has at least two distinct non-real
roots.

\end{teorem}

\nt Now we shall characterize graphs with few edge cover roots.
Note that zero is one of the roots of $E(G,x)$ with multiplicity
$\rho(G)$. The next theorem characterize all graph $G$ whose edge
cover polynomials have exactly one distinct root. Note that
$E(K_{1,n},x)= x^n$.

 \begin{teorem}{\rm(\cite{oboudi2})}
 Let $G$ be a graph. Then $E(G,x)$ has exactly one
distinct root if and only if every connected component of $G$ is
star.
\end{teorem}

\ms

\nt We need the following definition  to study graphs with two
distinct edge cover roots.

\nt  Let $H$ be a graph of order $n$ and size $m$. Suppose
$\{v_1,\ldots, v_n\}$ is the vertex set of $H$. By $H(r)$ we mean
the graph obtained by joining $r_i\geq 1$  pendant vertices to
vertex $v_i$, for $i = 1,\ldots n$ such that $\sum r_i= r$. If $m$
is the size of $H$, then $H(r)$ is a graph of order $n + r$ and
size $m + r$. The graph $C_4(9)$ shown in Figure \ref{figure4}.


\begin{teorem} \label{theorem20} {\rm(\cite{oboudi2})}
Let $G$ be a graph. Then $E(G,x)=x^r(x + 1)^m$, for some natural
numbers $r$ and $m$, if and only if there exists a graph $H$ with
size $m$ such that $G = H(r)$.

\end{teorem}

\nt The next theorem characterizes all graphs $G$ for which
$E(G,x)$ has exactly two distinct roots.

\begin{teorem} \label{theorem21}{\rm(\cite{oboudi2})}
 Let $G$ be a connected graph whose edge cover
polynomial has exactly two distinct roots. Then one of the
following holds:
\begin{enumerate}

\item[(i)] $G=H(r)$, for some connected graph $H$ and natural
number $r$.

\item[(ii)] $G = K_3$.

\item[(iii)]$\delta(G) = 1$, $E(G,x) = x^{\frac{m+s}{2}} (x +
2)^{\frac{m-s}{2}}$, where $s$ is the number of pendant vertices
of $G$.

\item[(iv)] $\delta(G)= 2$, $E(G,x)= x^{\frac{m}{2}}(x +
2)^{\frac{m}{2}}$ ,$ a_2(G) = \frac{m}{2}$, and $G$ has  cycle of
length $3$ or $5$.
\end{enumerate}
\end{teorem}

\nt By Theorems \ref{theorem20} and \ref{theorem21} we have the
following corollaries:

\begin{korolari}{\rm(\cite{oboudi1})}
Let $G$ be a connected graph. If $E(G,x)$ has exactly two distinct
roots, then $Z(E(G,x))=\{-1,0\}$, $Z(E(G,x))=\{-2,0\}$ or
$Z(E(G,x))=\{-3,0\}$.  Also $\delta(G)=1$ or $\delta(G)=2$.
\end{korolari}

\begin{korolari}{\rm(\cite{oboudi1})}
Let $G$ be a connected graph. If $Z(E(G,x))=Z(E(K_3,x))=\{-3,0\}$
if and only if $G=K_3$.
\end{korolari}

\nt We have the following corollary:

\begin{korolari}
Let $G$ be a graph with at most two distinct edge cover roots.
Then $Z(E(G,x))\subseteq \{-3,-2,-1,0\}$.
\end{korolari}

\nt There are infinite graphs $G$ with $\delta(G)=1$ and
$Z(E(G,x))\subseteq \{-2,-1,0\}$. To see some algorithms for
constructing of these kind of  graphs refer to \cite{oboudi1}. For
the case $\delta(G)=2$, we can see that if $G$ is a graph with
$\delta(G)=2$ and $Z(E(G,x))=\{-2,0\}$, then $G$ has cycle with
length $3$ or $5$ (see \cite{oboudi1}). For this case there is the
following conjecture:

\begin{konjek}{\rm(\cite{oboudi1})}
There isn't any graph $G$ with $\delta(G)=2$ and
$Z(E(G,x))=\{-2,0\}$.
\end{konjek}

\nt It is proven that if $G$ is a graph  without cycle of length
$3$ or $5$ and  $\delta(G)=2$, then $E(G,x)$ has at least three
distinct roots (see \cite{oboudi1}). Theorem \ref{theorem21} and
above conjecture implies that if $\delta(G)=2$ and $G\neq K_3$,
then $E(G,x)$ has at least three distinct roots.

\nt The following theorem gives the roots of edge cover polynomial
of paths and cycles:

\begin{teorem}{\rm(\cite{oboudi2})}
\begin{enumerate}
\item[(i)]
 For any integer $n$, $E(P_{n},x)$ has the following non-zero roots,
\[
-2\left(1+{\displaystyle \cos\frac{2s\pi}{n-1}}\right), \quad
s=1,2,\ldots,\lfloor \frac{n}{2}\rfloor-1.
\]

\item[(ii)]
 For any integer $n\geq3$, $E(C_{n},x)$ has the
following zeros,
\[
c_{s}^{(n)}=-2\left(1+{\displaystyle
\cos\frac{(2s+1)\pi}{n}}\right), \quad s=0,1,\ldots,\lfloor
\frac{n}{2}\rfloor-1.
\]
\end{enumerate}
\end{teorem}

\nt Also we have the following theorem for roots of $E(P_{n},x)$
and $E(C_{n},x)$.

\begin{teorem}{\rm(\cite{oboudi2})}
The  roots of the family $\{E(P_n,x)\}$ and $\{E(C_{n},x)\}$ are
real and dense in $(-4,0]$.
\end{teorem}

\nt {\bf Remark.} It is interesting that the non zeros  roots of
$I(C_n,x)$ are exactly the inverse of (non zeros) roots of
$E(C_n,x)$.

\nt Note that in \cite{oboudi2} proved that if every block of the
graph $G$ is $K_2$ or a cycle, then all real roots of $E(G, x)$
are in the interval $(-4, 0]$.

\ms

\nt For most of graphs polynomials such as chromatic polynomial,
matching polynomial, independence polynomial and characteristic
polynomial there is no constant bound for the roots (complex) of
these polynomials . Surprisingly,  in the following theorem we
observe that there is a constant bound for the roots of the edge
cover polynomials.

\begin{teorem}{\rm(\cite{oboudi2})}
All roots of the edge cover polynomial lie in the ball
$$\Big\{z\in \mathbb{C}: |z| < \frac{(2 + \sqrt{3})^2}{1
+\sqrt{3}}\simeq 5.099 \Big\}.$$
\end{teorem}

\section{Open problems and conjectures}

\nt In this section we state and review some open problems and
conjectures related to the subject of paper.

\begin{problem}{\rm(\cite{zero})}
 Characterize all graphs with exactly three distinct domination
roots $\{0,\frac{-3\pm \sqrt{3}i}{2}\}$.
\end{problem}

\begin{problem}{\rm(\cite{four})}
 Characterize all graphs with exactly four distinct domination
roots.
\end{problem}

\begin{problem}{\rm(\cite{zero})}
 Characterize all graphs with no real domination roots except
zero.
\end{problem}

\begin{konjek}{\rm(\cite{zero})}
 The set of integer domination roots of any graphs is a subset of
$\{-2,0\}$.
\end{konjek}

\begin{konjek}
Real roots of the families $D(P_{n},x)$ and $D(C_{n},x)$ are dense
in the interval $[-2,0]$, for $n\geq 4$.
\end{konjek}

\nt In this paper we obtained Theorem \ref{theorem18} for graphs
which its independence polynomials have few roots which is for
line graphs. But complete characterization remain as open problem:

\begin{problem}
 Characterize all graphs with few independence roots.
\end{problem}

\begin{konjek}{\rm(\cite{oboudi1})}
 If $\delta(G)=2$ and $G\neq K_3$, then $E(G,x)$ has at least
three distinct roots.
\end{konjek}

\begin{konjek}{\rm(\cite{oboudi1})}
 Let $G$ be a graph. Then $E(G,x)$ has at least $\delta(G)$
distinct roots.
\end{konjek}

\begin{konjek}{\rm(\cite{oboudi2})}
Let $G$ be a graph with $\delta(G) = 2$. If $E(G, x)$ has only
real roots, then all connected components of $G$ are cycles.
\end{konjek}

\end{document}